\documentstyle[12pt,amssymb]{article}
\newcommand{\llbracket}{[[}
\newcommand{\rrbracket}{]]}

\newtheorem{theo}{Theorem}
\newtheorem{prop}{Proposition}

\newcounter{remark}
\newenvironment{remark}{\stepcounter{remark}\noindent{\bf Remark \arabic{remark}}}{\vskip2ex}
\newcommand{\qed}{\hfill $\Box$\vskip2ex}

\newcommand{\bR}{{\Bbb R}}

\newcommand{\frakh}{{\frak h}}
\newcommand{\frakg}{{\frak g}}
\newcommand{\frakk}{{\frak k}}
\newcommand{\frakp}{{\frak p}}
\newcommand{\frakr}{{\frak r}}

\newcommand{\Ad}{{\rm Ad}}

\newcommand{\supp}{{\rm supp}}
\newcommand{\levy}{L\'{e}vy }
\newcommand{\Dom}{{\rm Dom}}

\begin{document}

\begin{center}
{\Large\bf Isotropic self-similar Markov processes}
\vspace{2ex}

 {\large Ming Liao}\footnote{Department of Mathematics, Auburn
   University, Auburn, AL 36849, USA. Email: liaomin@auburn.edu}
 {\large and Longmin Wang}\footnote{School of Mathematical Sciences, Nankai
    University, Tianjin, China. Email: wanglm@nankai.edu.cn

    \ \ The second author was partially supported by NSF of China (No.
    10871103).}
\end{center}

\begin{quote}
{\bf Summary } We show that an isotropic self-similar Markov process
in $\bR^d$ has a skew product structure if and only if its radial
and angular parts do not jump at the same time.


\end{quote}

\section{Introduction and main result}\label{sec:intro}

It is well known that a Brownian motion in $\bR^d$ ($d\ge2$)
has a skew product structure, that is, it may be expressed as a product
of its radial process and a time changed spherical Brownian motion. Moreover,
the radial process is a Bessel process and is independent of
the spherical Brownian motion, and the time change is adapted
to the radial process. This decomposition is naturally related to
the invariance of the Brownian motion under the group $O(d)$
of orthogonal transformations on $\bR^d$. More generally, Galmarino~\cite{Galmarino1962}
proved that a continuous isotropic or $O(d)$-invariant Markov process in
$\bR^d$ is also a skew product of its radial motion and
an independent spherical Brownian motion with
a time change. Pauwels and Rogers~\cite{PauwelsRogers1988} and
Liao~\cite{Liao2009a} extended these results to more general settings.

Because any continuous isotropic Markov process has a skew product
  structure, it is therefore natural to consider a similar skew
product for discontinuous isotropic Markov processes. Graversen and
Vuolle-Apiala~\cite{GraversenVA1986} discussed a skew product for
isotropic $\alpha$-self-similar Markov processes, which include the
purely discontinuous symmetric $(1/\alpha)$-stable processes. Their
main result says that after a time change due to
Lamperti~\cite{Lamperti1972} and Kiu~\cite{Kiu1980}, the radial
process and the angular process are respectively multiplicatively
invariant and $O(d)$-invariant Markov processes, and are
independent. This leads to a skew product structure similar to that
of a Brownian motion. However, as will be shown later, the
independence part of this interesting result holds only under a
rather restrictive condition, which excludes for example the
  symmetric $(1/\alpha)$-stable processes for $\alpha>1/2$.
 We note that the proof of Proposition~2.4 in \cite{GraversenVA1986} has an
error in the conditional expectation argument.

The aim of this paper is to clarify this rather important point. We
will show that an isotropic $\alpha$-self-similar Markov process has
a skew product structure if and only if its radial and angular parts
  do not jump at the same time.

After the completion of this paper, we found the independence error
in \cite{GraversenVA1986} has been noted in Bertoin and Werner
\cite{bw}. They used the correct part of \cite{GraversenVA1986} in
their work, but did not pursue the independence problem.

We now describe our setup more precisely. All processes considered
in this paper are assumed to have c\`adl\`ag paths (right continuous
paths with left limits). Let $x_t$ be a (time homogeneous) Markov
process in $\bR^d$, $d\geq 2$, with transition function $P_t$
 satisfying the usual simple Markov property. We will allow process $x_t$ to have a possibly finite life time,
and as usual let $P_x$ denote the distribution of process starting
at $x$ on the canonical space of c\`adl\`ag paths with possibly
finite life times.

The restriction of process $x_t$ on $E=\bR^d-\{0\}$, defined before reaching
the hitting time of origin $0$, is also a Markov process. For simplicity,
this process together with its transition function and distribution are still denoted
by $x_t$, $P_t$ and $P_x$.

The process $x_t$ in $\bR^d$ is said to be \emph{isotropic} or \emph{$O(d)$-invariant}
if
\begin{equation}
  P_t(\phi(x),\phi(B))=P_t(x,B)
  \label{eq:isotropic}
\end{equation}
for any $\phi\in O(d)$, $x\in\bR^d$ and Borel subset $B\subset\bR^d$.
This is equivalent to saying that process $\phi(x_t)$ with $x_0=x$ has
the same distribution as process $x_t$ with $x_0=\phi(x)$.

The process $x_t$ is said to be \emph{$\alpha$-self-similar}, or \emph{$\alpha$-s.s.} in short,
for some constant $\alpha>0$, if
\begin{equation}
  P_{\lambda t}(x,B)=P_t (\lambda^{-\alpha} x,\lambda^{-\alpha} B)
  \label{eq:selfsimilar}
\end{equation}
for any $\lambda>0$, $x\in\bR^d$
and $B\subset\bR^d$. This is equivalent to saying
that process $x_{\lambda t}$ with $x_0= x$ has the same distribution
as process $\lambda^\alpha x_t$ with $x_0=\lambda^{-\alpha}x$.

It is clear that if $x_t$ is $O(d)$-invariant and/or
$\alpha$-self-similar, so is its restriction to $E$.

In what follows, we will exclusively consider an isotropic Markov
process $x_t$ in $E$. For $x=(x_1,\ldots,x_d)\in\bR^d$, we denote by
$\vert x\vert=\sqrt{x_1^2+\cdots+x_d^2}$ the radial part of $x$. Let
$r_t=\vert x_t\vert$ and $\theta_t=x_t/r_t$ be the radial and
angular parts of $x_t$. It is easy to show using the
$O(d)$-invariance (see for example \cite{Liao2009a}) that $r_t>0$ is
a $1$-dim Markov process with transition function $R_t$ given by
\[R_tf(r) = P_t(f\circ\pi_1)(x)\]
for $r>0$ and Borel function $f$ on $(0,\,\infty)$, where $\pi_1$: $E
\to(0,\,\infty)$ is the natural projection given by $y\mapsto\vert y\vert$
and $x$ is any point of $E$ with $r=\vert x\vert$.
As the angular process $\theta_t$
lives in the unit sphere $S^{d-1}$ which is invariant under the action
of $O(d)$, one would expect that it should inherit the $O(d)$-invariance
of $x_t$ in some sense. This leads to the following definition
of a skew product structure.
\vspace{2ex}

\noindent {\bf Definition} \  Let $x_t$ be an isotropic Markov process
in $E$. We say that $x_t$ has a \emph{skew product structure}
if $x_t=r_t\xi_{A_t}$,
where $A_t$ is a continuous and strictly increasing process with $A_0=0$,
adapted to the radial process $r_t$,
and $\xi_t$ is an $O(d)$-invariant Markov process in $S^{d-1}$
and is independent of process $r_t$.
\vspace{2ex}

Because $O(d)$ acts on $S^{d-1}$ transitively, $S^{d-1}$ may
be regarded as a homogeneous space of $O(d)$.
Invariant Markov processes in homogeneous spaces are Feller processes,
and their generators may be expressed explicitly in terms
of an invariant differential operator and a \levy measure (see
Section~\ref{sec:proof} for more details), thus providing
a useful tool for studying these processes.

The following is our main theorem.
\begin{theo}\label{theo:main}
  Let $x_t$ be an isotropic $\alpha$-self-similar Markov process
  in $E=\bR^d-\{0\}$ ($d\ge 2$). Then $x_t$ has a skew product structure if and only if
  its radial and angular parts do not jump at the same time, that is,
  for all $x\in E$, $P_x$-almost surely, $r_t$ and $\theta_t$ cannot
  jump together at same time $t$ for any $t\geq 0$.
\end{theo}

\noindent {\bf Proof of the necessity part\ }
  The sufficiency of the condition will be proved in
  Section~\ref{sec:proof}.
  For the necessity, we assume that $x_t$ has the skew product $x_t = r_t  \xi_{A_t}$.
  Since $r_t$ is c\`adl\`ag, by
  \cite[Proposition~I.1.32]{JacodShiryaev1987}, the random set
  $\{\Delta r_t\not=0\}$ is thin in the sense that there is a
  sequence of stopping times $\tau_n$ such that $\{\Delta
  r_t\not=0\} = \bigcup_n \llbracket\tau_n\rrbracket$, where $\Delta
  r_t = r_t - r_{t-}$, and
  $\llbracket\tau_n\rrbracket$ is the graph of $\tau_n$, i.e.,
  $\llbracket\tau_n\rrbracket = \{(\omega,t),\ t\in\bR_+,\
  t=\tau_n(\omega)\}$.  For any $n\ge1$, the time
  $A_{\tau_n}$ is measurable in process $r_t$, and the
  independence of $r_t$ and $\xi_t$ implies that
  $A_{\tau_n}$ is independent of $\xi_t$. As a Feller process, $\xi_t$
  is quasi-left-continuous. In particular, $\xi_t$ does
  not jump at a fixed time, and it is easy to see that $\xi_t$ does not
  jump at $A_{\tau_n}$. This implies that the radial part $r_t$
  and the angular part $\theta_t=\xi_{A_t}$ of $x_t$ do not jump
  simultaneously.
  \qed

\begin{remark}
  Note that an isotropic self-similar Markov process
  may not satisfy the condition in Theorem~\ref{theo:main}. The most
  famous examples are the symmetric $(1/\alpha)$-stable \levy processes
  for $\alpha>1/2$. Their \levy measures are absolutely continuous on
  $\bR^d-\{0\}$, so their radial and angular parts may jump together,
  and thus do not possess a skew product structure as defined above.
  On the other hand, we will see later that there are many isotropic
  $\alpha$-s.s. Markov processes that do possess a skew product structure.
\end{remark}

\begin{remark}
  It is evident that the proof above is also valid for
  a general isotropic Markov process. That is, the jump condition in
  Theorem~\ref{theo:main} is also necessary for an isotropic
  Markov process to have a skew product structure.
\end{remark}

The rest of this paper is devoted to proving the sufficiency part
of Theorem~\ref{theo:main}. In Section~\ref{sec:timechange},
we will recall the time change used in \cite{GraversenVA1986}, and we
will show that $x_t$ is $\alpha$-s.s..
if and only if the time changed process is invariant
under the scalar multiplication.
The key fact is that if $x_t$ is $\alpha$-s.s., then the time changed process is invariant under
a transitive group on $E$ and hence may be viewed as an invariant Markov
process in a homogeneous space. Under this viewpoint, we complete
the proof of our main theorem in Section~\ref{sec:proof}.

\section{Time changed processes}\label{sec:timechange}

Let $x_t$ be an isotropic Markov process (not necessarily
$\alpha$-s.s.) in $E=\bR^d-\{0\}$ with a possibly finite life time
$\xi$. Fix $\alpha>0$. The following random time change was
introduced by Lamperti~\cite{Lamperti1972} for $\bR_+$-valued
processes. Define
\begin{equation}
  \label{eq:CAF}
  A_t = \int_0^t |x_{s}|^{-1/\alpha}ds,
\end{equation}
which is a continuous and strictly increasing function for $t<\xi$.
Its inverse $T_t$ is given by
\begin{equation}
  \label{eq:timechange}
  T_t = \inf\{s\ge0; \ A_s\ge t\},\ t< A_{\xi-}.
\end{equation}

We define a new process $\{\bar x_{t}\}$ by $\bar x_{t} = x_{T_t}$
  for $t<A_{\xi-}$ and $x_t=\Delta$ otherwise, where $\Delta$ is a cemetery
point added to $E$. By Theorem~10.11 of~\cite{Dynkin1965}, $\bar
x_{t}$ is also a time homogeneous Markov process with c\`adl\`ag
paths. Let $\bar P_t(x,B)$ be the transition function of $\bar
x_{t}$. Note that $\bar x_t$ and $\bar P_t$ are also isotropic.

It is easy to show that
\begin{equation}\label{eq:Tt}
  T_t = \int_0^t |\bar x_u|^{1/\alpha} du
\end{equation}
for $t<A_{\xi-}$. That is, $T_t$ is determined by the time changed
process $\bar x_t$ and is also continuous and strictly increasing.
Note that $A_t$ is the inverse of $T_t$. Thus we may start with an
isotropic Markov process $\bar x_t$ in $E$ and recover the original
process $x_t$ as $\bar{x}_{A_t}$.

The process $\bar x_t$ is said to be \emph{multiplicatively invariant}, if
\begin{equation}\label{eq:m_inv}
  \bar P_t(x,B) = \bar P_t(\lambda x,\lambda B)
\end{equation}
for any $\lambda >0$, $x\in E$ and Borel subset $B\subset E$. This is
equivalent to saying that process $\lambda \bar x_t$ with $\bar x_0=x$
has the same distribution as process $\bar x_t$ with $\bar x_0 =
\lambda x$.

The following theorem relates the $\alpha$-self-similarity of $x_t$ to
the multiplicative invariance of $\bar x_t$. The multiplicative
invariance of $\bar x_t$ was proved by Kiu~\cite{Kiu1980}, but the
present proof is simpler and more probabilistic, and also establishes
its converse.
\begin{theo}\label{theo:kiu}
  The process $x_t$ is $\alpha$-s.s. if and only if the time changed
  process $\bar x_t$ is multiplicatively invariant.
\end{theo}
\noindent {\bf Proof \ }
For simplicity, we will work on the canonical probability space of
c\`adl\`ag paths with possibly finite life time. We will also  write
$x_\cdot$ for a path $x_t$ in $E$ and $x_{\lambda\cdot}$ for path
$t\mapsto x_{\lambda t}$ for $\lambda>0$. To indicate the dependence
on a path $x_\cdot$, we will write $A_t(x_{\cdot})$ and
$T_t(x_{\cdot})$ instead of $A_t$ and $T_t$.

Assume that $x_t$ is $\alpha$-s.s.. Then $(\lambda^\alpha
x_{\lambda^{-1} t}, P_{\lambda^{-\alpha}x})$ is the same Markov
process as $(x_t,P_x)$ and consequently, under $P_{\lambda^{-\alpha} x}$, the
distribution of $(\lambda^{\alpha} x_{\lambda^{-1}t},
T_t(\lambda^{\alpha}x_{\lambda^{-1}\cdot}))$ equals that of
$(x_t,T_t(x_\cdot))$ under $P_x$.
Since
\[ A_t(\lambda^\alpha x_{\lambda^{-1}\cdot}) = \lambda^{-1}\int_0^t
|x_{\lambda^{-1}s}|^{-1/\alpha} ds = \int_0^{\lambda^{-1} t}
|x_{s}|^{-1/\alpha} ds = A_{\lambda^{-1}t}(x_{\cdot}), \]
we obtain that $T_t(\lambda^\alpha x_{\lambda^{-1}\cdot}) = \lambda
T_t(x_\cdot)$. Note that the processes $\lambda^\alpha \bar x_t$ and
$\bar x_t$ are respectively measurable functionals of the processes
$(\lambda^\alpha x_{\lambda^{-1}t},\lambda T_t(x_\cdot))$ and
$(x_t,T_t(x_\cdot))$ of the same form.
It follows that process $\lambda^{\alpha}\bar x_t$ with
$\bar x_0 = \lambda^{-\alpha}x$ has the same distribution as process
$\bar x_t$ with $\bar x_0 = x$. This proves the multiplicative
invariance of $\bar x_t$.

Conversely, assume that $\bar x_t$ is multiplicatively invariant. Then
the process $\lambda^\alpha \bar x_t$ with $\bar x_0 = x$ has the same
distribution as the process $\bar x_t$ with $\bar x_0 = \lambda^\alpha
x$.
Let $\bar{T}_t(\bar{x}_\cdot)$ denote the integral in (5) and
let $\bar{A}_t(\bar{x}_\cdot)$ be its inverse as a function of $t$.
Then $A_t(x_\cdot)=\bar{A}_t(\bar{x}_\cdot)$, and the distribution
of $(\lambda^\alpha \bar{x}_t,\bar{A}_t(\lambda^\alpha \bar{x}_\cdot))$
with $\bar{x}_0=x$ equals that of $(\bar{x}_t,\bar{A}_t(\bar{x}_\cdot))$
with $\bar{x}_0=\lambda^\alpha x$.  Because
$\bar{T}_t(\lambda^\alpha\bar{x}_\cdot)=\lambda\bar{T}_t(\bar{x}_\cdot)$,
$\bar{A}_t(\lambda^\alpha\bar{x}_\cdot)=\bar{A}_{\lambda^{-1}t}(\bar{x}_\cdot)
= A_{\lambda^{-1}t}(x_\cdot)$. The $\alpha$-self-similarity of $x_t$
now follows from a substitution of $\bar{A}_t(\lambda^\alpha
\bar{x}_\cdot)$ for $t$ in $\lambda^\alpha\bar{x}_t$.
\qed

As in~\cite{GraversenVA1986}, the semigroup property implies that
there is a $\gamma\ge0$ such that $\bar P_t(x,E) = e^{-\gamma t}$ for
$t\ge0$ and $x\in E$. When $\gamma>0$, $\bar x_t$ will have a finite
life time, or equivalently, $\bar P_t$ is not
conservative. But we may define a new transition function $\hat P_t$
by
\[  \hat P_t(x,B) = e^{\gamma t} \bar P_t(x,B), \quad t\ge0,\ x\in E,\
  B\subset E. \]
Then $\hat{P}_t$ is a conservative transition function, and
the associated conservative Markov process $\hat{x}_t$ is isotropic
and multiplicatively invariant.  The process $\bar{x}_t$ is just
process $\hat{x}_t$ killed at an independent exponential time
of rate $\gamma$.

\section{Proof of the sufficiency part in Theorem~\ref{theo:main}}\label{sec:proof}

Let $d\ge2$ and let $GL(d,\bR)$ be the group of the nonsingular linear
transformations on $\bR^d$. Let $G$ be the similarity group of
$\bR^d$, that is,
\[ G = \{ g\in GL(d,\bR); \ |gv| = |g| |v| {\rm\ for\ any\ } v\in
\bR^d \}, \] where $|v| = \sqrt{v_1^2+\cdots +v_d^2}$ for
$v=(v_1,\ldots,v_d)\in \bR^d$ and $|g|$ is the operator norm of
$g\in GL(d,\bR)$, that is, $|g| = \sup_{|v|=1}|gv|$. For $c>0$,
define the linear transformation $m_c$ by $m_cv=cv$ for $v\in\bR^d$.
Let $R=\{m_c;\ c>0\}$ and $H=O(d)$. Then $R$ and $H$ are both normal
subgroups of $G$. Moreover, $G$ is the direct product of $R$ and
$H$.

Note that $G$ acts transitively on $E=\bR^d-\{0\}$. Fix
$o=(0,\ldots,0,1)$. The subgroup of $G$ fixing $o$ is $K=O(d-1)$.
We may identify $G/K$ with $E$ via the map $gK \mapsto go$, $H/K$
with the sphere $S^{d-1}$ via $hK \mapsto ho$, and $R$ with a ray
in $E\subset \bR^d$ via $r \mapsto ro$. Note that $E$ is diffeomorphic to
the product space $R\times S^{d-1}$.

The reader is referred to section 2.2 of~\cite{Liao2004} for some basic definitions
about invariant Markov processes in homogeneous spaces.
Let $\frakg$, $\frakr$,
$\frakh$ and $\frakk$ be respectively the Lie algebras of $G$, $R$, $H$ and
$K$. There is an $\Ad(K)$-invariant subspace $\frakp$ such
that $\frakh = \frakk \oplus\frakp$.
Then the exponential map of $G$ provides a natural
local diffeomorphism from $\frakr \oplus \frakp$ to $E$. Let $n=\dim G$ and
let $\{X_1,\ldots,X_n\}$ be a basis of $\frakg$ such that
$X_1\in \frakr$, $X_2,\ldots,X_d\in\frakp$ and
$X_{d+1},\ldots,X_n\in\frakk$. Let $\pi$: $G \to E$ be the map $g
\mapsto go$. Restricted to a sufficient small
neighborhood $V$ of $0$, the map
\[ \phi:\ \bR^d \ni y=(y_1,\ldots,y_d) \mapsto \pi(e^{\sum_{j=1}^d
  y_jX_j}) \in E\]
is a diffeomorphism and $y_1,\ldots,y_d$ may be used as local
coordinates on $\phi(V)$. As in Section 2.2 of Liao~\cite{Liao2004},
we may extend $y_j$ to $E$ such that $y_j\in C_c^\infty(E)$ (the
space of smooth functions on $E$ with compact supports) and for any
$x\in E$, $k\in K$,
\begin{equation}
  \label{eq:localcrd}
  \sum_{j=1}^d y_j(x) \Ad(k)X_j = \sum_{j=1}^d y_j(kx) X_j.
\end{equation}

As in Section~\ref{sec:timechange}, we let $x_t$ be an isotropic
$\alpha$-s.s. Markov process starting at $x\in E$. Recall that the
time changed process $\bar x_t$ defined before is a $G$-invariant
Markov process in $E$ with transition function $\bar P_t$ (see
Theorem~\ref{theo:kiu}). Thus for any $f\in C_c^\infty(E)$ and $x\in
E$,
\[  \bar P_tf(x) = \bar P_t(f\circ g)(o), \]
where $g\in G$ is chosen to satisfy $x=go$. As an easy consequence,
$\bar P_t$ is a $G$-invariant Feller semigroup on $E$.

Let $L$ be the generator of $\bar x_t$ with domain $\Dom(L)$. An
explicit formula for the generator of an invariant Markov process in
a homogeneous space was obtained by Hunt \cite{hunt}. By Theorem~2.1
of \cite{Liao2004}, which is a more convenient version of Hunt's
formula, $\Dom(L)$ contains $C_c^\infty(E)$ and for $f\in C_c^\infty
(E)$,
\begin{equation}
  \label{eq:generator}
  Lf(o) = Tf(o) + \int_E \left[ f(x) - f(o) - \sum_{j=1}^d y_j(x)
    \frac{\partial }{\partial y_j} f(o)\right] \Pi(dx),
\end{equation}
where $T$ is a $G$-invariant diffusion generator and $\Pi$ is a
$K$-invariant \levy measure on $E$. There exist a $d\times d$
non-negative definite symmetric matrix $(a_{ij})$ and constants $c_i$
such that for $f\in C_c^\infty(E)$,
\begin{equation}
  \label{eq:diffusion}
  Tf(o) = \frac{1}{2} \sum_{i,j=1}^da_{ij}X_i^lX_j^l(f\circ\pi)(e) +
  \sum_{i=1}^d c_i X_i^l(f\circ\pi)(e),
\end{equation}
where $X_i^l$ is the left invariant vector field
on $G$ determined by $X_i$.
Moreover, the coefficients $a_{ij}$ and $c_i$ satisfy
\begin{equation}
  \label{eq:diffcoef}
  a_{ij} = \sum_{p,q=1}^d a_{pq}b_{ip}(k)b_{jq}(k) {\rm\quad and\quad } c_i =
  \sum_{p=1}^d c_p b_{ip}(k), \quad \forall k\in K,
\end{equation}
where the orthogonal matrix $(b_{ij}(k))$ is determined by $\Ad(k)X_j
= \sum_{i=1}^d b_{ij}(k) X_i$ for $j=1,\ldots,d$.

  Since $R$ commutes with $H$ and $\frakp$ is $\Ad(K)$-invariant,
$\Ad(k)X_1=X_1$ and
$\Ad(k)X_i\in\frakp$ for $i\ge2$ and $k\in K$. Thus $b_{11}(k)=1$,
$b_{i1}(k)=b_{1i}(k)=0$ for $i\ge2$. Then (\ref{eq:diffcoef}) implies
that
\[ a_{i1} = \sum_{p=2}^d a_{p1} b_{ip}(k), \ a_{1i} = \sum_{q=2}^d
a_{1q} b_{iq}(k), {\rm\quad and\quad } c_i = \sum_{p=2}^d c_p
b_{ip}(k), \quad 2\le i\le d.  \] In other words, the vectors
$X=\sum_{i=2}^d a_{i1}X_i$ and $Y=\sum_{i=2}^d a_{1i}X_i$ are
invariant under the action of $\Ad(k)$ for all $k\in K$, which
implies that $X=Y=0$. Hence $a_{i1}=a_{1i}=0$ for $2\le i\le d$. The
operator $T_2$ defined by
\[ T_2f(o) = \frac{1}{2} \sum_{i,j=2}^d a_{ij}X_i^lX_j^l(f\circ
\pi)(e) + \sum_{i=2}^d c_i X_i^l(f\circ \pi)(e),\ f\in C_c^\infty(E)
\] may be viewed as an $H$-invariant diffusion generator on the sphere
$S^{d-1}=H/K$. It is well known that there is a constant $c\ge0$
such that $T_2=c\Delta$, where $\Delta$ is the Laplace-Beltrami
operator on $S^{d-1}$. We define the diffusion generator $T_1$ by
\[ T_1f(o) = \frac{1}{2} a_{11}X_1^lX_1^l(f\circ \pi)(e) +
c_1X_1^l(f\circ \pi)(e),\ f\in C_c^\infty(E). \]
Note that operator $T_1$ acts along $R$.
We have proved that $T=T_1+T_2$ in the sense that
\[ Tf(r,\theta) = (T_1f(\cdot,\theta))(r) + (T_2f(r,\cdot))(\theta) \]
for $r\in R$ and $\theta \in S^{d-1}$.

Let $\pi_1$ (resp. $\pi_2$) be the
projection from $E$ onto $R$ (resp. $S^{d-1}$). Then for $x\in E$,
$\pi_1(x)$ (resp. $\pi_2(x)$) may be identified with $|x|$ (resp. $
x/|x|$).
Let $\rho_t=\pi_1(\bar x_t)$ and $\xi_t=\pi_2(\bar x_t)$. By the $O(d)$-invariance of
$\bar x_t$, $\rho_t$ is a \levy process on $R$ starting at $\pi_1(x)$ and $\xi_t$ is an
$O(d)$-invariant Feller process on the sphere $S^{d-1}$ starting at $\pi_2(x)$.
\begin{prop}\label{prop:levymeasure}
  $\rho_t$ and $\xi_t$ are independent if and only if the \levy measure
  $\Pi$ of $\bar x_t$ is concentrated on $R\bigcup S^{d-1}$, where
  $R$ and $S^{d-1}$ are regarded as subsets of $E\cong R\times
  S^{d-1}$.
\end{prop}
\noindent{\bf Proof }
  Assume that $\rho_t$ and $\xi_t$ are independent. Let $f_1$ (resp.
  $f_2$) be a smooth function on $R$ (resp. $S^{d-1}$) vanishing near $o$.
  Let $f(x) = f_1(\pi_1(x))f_2(\pi_2(x))$. Then by
  (\ref{eq:generator}), $\Pi(f)=\int_{E}f(x)\Pi(dx) = Lf(o)$. From the
  independence of $\rho_t$ and $\xi_t$, we have that
  \[ Lf(o) = \lim_{t\to0} \frac{E[f(\bar x_t)]}{t} = \lim_{t\to0}
  \frac{E[f_1(\rho_t)] E[f_2(\xi_t)]}{t}. \]
  It follows that $\Pi(f)=0$ since $E[f_1(\rho_t)] = tO(t)$ and
  $E[f_2(\xi_t)] = tO(t)$ as $t\to0$.

  Now fix a point $x=(r,\theta)\in E$ such that $r$ and
  $\theta$ are not the point $o$. We may choose positive functions
  $f_1$ on $R$ and $f_2$ on $S^{d-1}$ satisfying the above conditions and
  additionally, we assume that $f_1=1$ near $r$ and that $f_2=1$ near
  $\theta$. Then there exists a neighborhood $U$ of $(r,\theta)$ such that
  $\Pi(U) \le \Pi(f_1f_2) = 0$. Hence $(r,\theta)$ is not contained in
  the support of $\Pi$. It follows that $\supp \Pi \subset R\bigcup
  S^{d-1}$.

  Conversely, let $\Pi=\Pi_1 + \Pi_2$ be such that $\Pi_1$ and $\Pi_2$
  are respectively \levy measures on $R$ and $S^{d-1}$, regarded as
  measures on $E$ supported by $R$ and $S^{d-1}$. For $i=1$, $2$,
  let $L_i$ be generators with diffusion parts $T_i$ and \levy measures
  $\Pi_i$.
  Our computation shows that $L=L_1+L_2$ at point $o$, and by
  the $G$-invariance of the three operators, $L=L_1+L_2$ on $E$.
  Note that when restricted to $R$ (resp. $S^{d-1}$), $L_1$
  (resp. $L_2$) is the generator of $\rho_t$ (resp. $\xi_t$).
Let $\tilde{\rho}_t$ be a \levy process in $R$ with generator $L_1$
and let $\tilde{\xi}_t$ be an $O(d)$-invariant Markov process in
$S^{d-1}$ with generator $L_2$, and let them be independent. Then
$\tilde{x}_t=(\tilde{\rho}_t,\tilde{\xi}_t)$ is a $G$-invariant
Markov process in $E$ with generator $L=L_1+L_2$. By the uniqueness
in Theorem~2.1
 of \cite{Liao2004}, the processes $\tilde{x}_t=(\tilde{\rho}_t,\tilde{\xi}_t)$
and $\bar{x}_t=(\rho_t,\xi_t)$ have the same distribution because
they have the same generator. This shows that $\rho_t$ and $\xi_t$
are independent.
\qed

\noindent{\bf Proof of sufficiency in Theorem~\ref{theo:main}}
  Now we assume that the radial and angular parts of $x_t$ do not jump at
  same time.
  It is obvious that the time change $\bar x_t
  = x_{T_t}$ does not change the directions of jumps. Thus the \levy measure of
  $\bar x_t$ is concentrated on the radial and angular axes.
  By Proposition~\ref{prop:levymeasure}, $\rho_t=\pi_1(\bar x_t)$ and
  $\xi_t = \pi_2(\bar x_t)$ are independent, and $\bar x_t = \rho_t\xi_t$.
  Recall that $A_t = \int_0^t |x_s|^{-1/\alpha}ds$ is the
  inverse of $T_t = \int_0^t (\rho_s)^{1/\alpha} ds$ and $x_t=\bar{x}_{A_t}$. Then $x_t = r_t
  \xi_{A_t}$, where $r_t = |x_t| = \rho_{A_t}$ is an
  $\alpha$-self-similar process on $(0,\ \infty)$ and is independent of
  $\xi_t$. Thus the skew product structure of $x_t$ is established.
\qed

\begin{remark}
  Our proof shows that the time change by $A_t$ provides a 1-1
  correspondence between isotropic $\alpha$-s.s. Markov processes and
  $G$-invariant Markov processes in $E$. Thus, given any \levy measure
  supported by $R\bigcup S^{d-1}$, there is a unique isotropic $\alpha$-s.s.
  Markov process in $E$ that possesses a skew product structure.
  We also note that any isotropic $\alpha$-s.s. Markov process has the
  strong Markov property, because the strong Markov property is
  possessed by the time changed process and is preserved by the
  inverse time change.
\end{remark}

\bibliographystyle{abbrv}

\end{document}